\documentclass[12pt]{article}
\usepackage[mathscr]{eucal}
\usepackage{amsmath,amssymb,amscd}
\usepackage{color}
\usepackage{epsfig}

\def \beq{ \begin{equation}}
\def \eeq{\end{equation}}

\begin{document}

\begin{center}
{\bf A proof of the butterfly theorem using the scale factor between
the two wings.}
\end{center}

\begin{center}
{\bf Martin Celli}
\end{center}

\begin{center}
October 21st 2016
\end{center}

\begin{center}
Departamento de Matem\'aticas\\
Universidad Aut\'onoma Metropolitana-Iztapalapa\\
Av. San Rafael Atlixco, 186. Col. Vicentina. Del. Iztapalapa. CP
09340. M\'exico, D.F.\\
E-mail: cell@xanum.uam.mx
\end{center}

\vskip1cm

{\bf Abstract.} We give a new proof of the butterfly theorem, based
on the use of several expressions involving the scale factor between
the two wings.\\

The aim of this article is to propose a new proof of the following
theorem:\\

{\bf Butterfly theorem.} {\it  Let \(M\) be the midpoint of a chord
\(PQ\) of a circle, through which two other chords \(AB\) and \(CD\)
are drawn. Let us assume that \(A\) and \(D\) do not belong to a
same half-plane defined by \(PQ\). Let \(X\) (respectively \(Y\)) be
the intersection of \(AD\) (respectively \(BC\)) and \(PQ\). Then
\(M\) is also the midpoint of \(XY\).}\\

\begin{figure}[h]
 \centering
    \includegraphics[width=15cm]{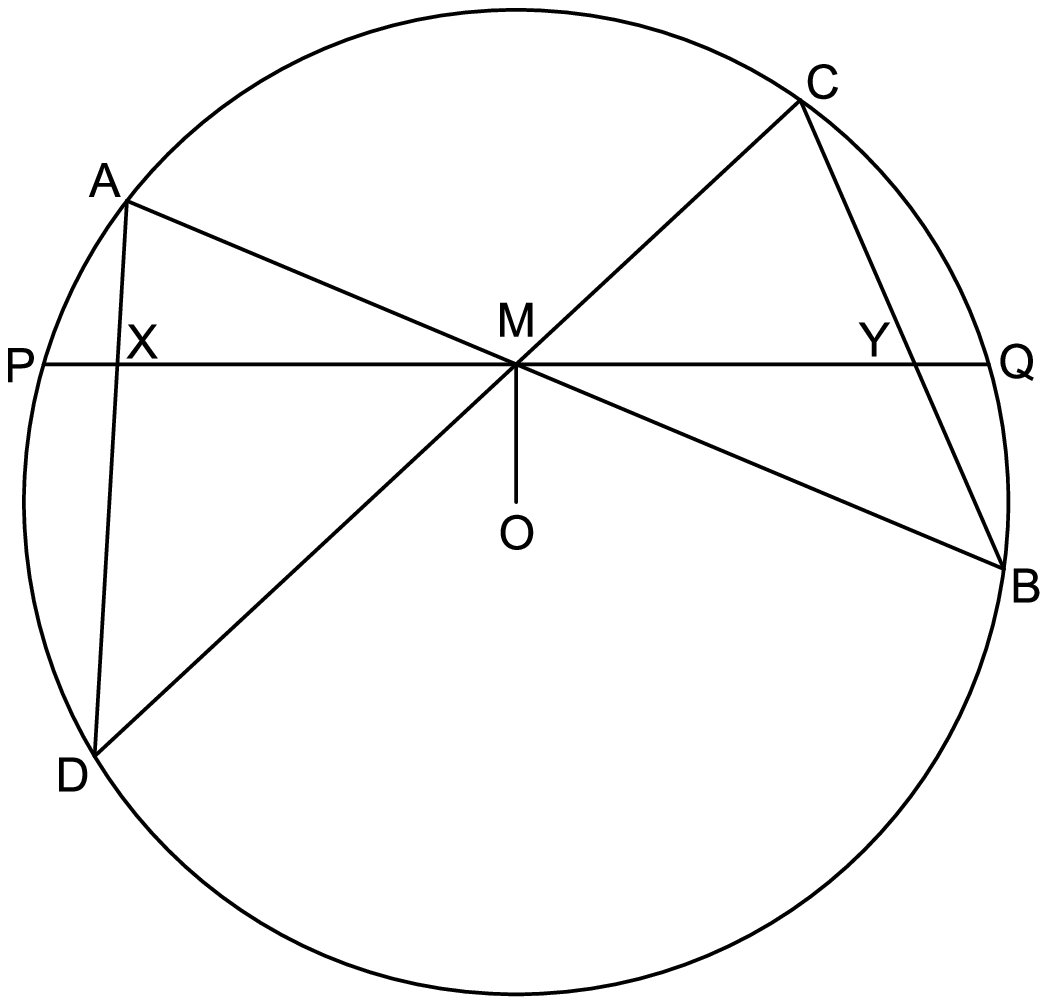}
      \end{figure}

Let \(O\) be the center of the circle. The points \(A\) and \(C\)
belong to a same half-plane defined by \(PQ\), the points \(B\) and
\(D\) to the other half-plane. For sake of simplicity, we can assume
that \(O\) belongs to the same half-plane as \(B\) and \(D\).\\

Several classic and recent proofs of this theorem are known ([1],
[2]). In our proof, we directly show that the ratio \(\sin
(CYM)/\sin (AXM)\) is nothing but the scale factor between the two
wings \(AMD\) and \(CMB\), which are similar by the inscribed angle
theorem:
\[\frac{\sin (CYM)}{\sin (AXM)}=\alpha
\mbox{, where }\alpha
=\frac{CM}{AM}=\frac{BM}{DM}=\frac{CB}{AD}\cdot \] More precisely,
we have:
\[\sin (AXM)=\frac{DM^2-AM^2}{AD.OM}\cdot \]
As a matter of fact:
\[DM^2-AM^2=||\overrightarrow{DO}+\overrightarrow{OM}||^2
-||\overrightarrow{AO}+\overrightarrow{OM}||^2\]
\[=OD^2+OM^2+2\overrightarrow{DO}.\overrightarrow{OM}
-(OA^2+OM^2+2\overrightarrow{AO}.\overrightarrow{OM})\]
\[=2\overrightarrow{DA}.\overrightarrow{OM} \mbox{ as } OA=OD\]
\[=2AD.OM\sin(AXM) \mbox{ as } OMX=OMY=\frac{OMX+OMY}{2}=\frac{\pi}{2},\]
because triangles \(OMP\) and \(OMQ\) are congruent. Similarly, we
have:
\[\sin (CYM)=\frac{BM^2-CM^2}{CB.OM}\cdot \]
Thus:
\[\frac{XM}{YM}=\frac{AM}{CM}\times \frac{CM}{YM}\times \frac{XM}{AM}
=\frac{AM}{CM}\times \frac{\sin (CYM)}{\sin (YCM)}\times \frac{\sin
(XAM)}{\sin (AXM)}\] by the law of sines, applied to triangles
\(AXM\) and \(CYM\)
\[=\frac{AM}{CM}\times \frac{\sin (CYM)}{\sin (AXM)}
\mbox{ by the inscribed angle theorem}\]
\[=\frac{AM}{CM}\times \frac{BM^2-CM^2}{CB.OM}\times \frac{AD.OM}{DM^2-AM^2}
=\frac{AM}{CM}\times \frac{AD}{CB}\times
\frac{BM^2-CM^2}{DM^2-AM^2}\]
\[=\frac{1}{\alpha}\times
\frac{1}{\alpha}\times \alpha ^2=1\cdot \]

\(\)\\

{\bf References.}\\
\(\mbox{[1]}\) Alexander Bogomolny, Butterfly theorem. Interactive
Mathematics Miscellany and Puzzles:\\
http://www.cut-the-knot.org/pythagoras/Butterfly.shtml\\
\(\mbox{[2]}\) Cesare Donolato, A proof of the butterfly theorem
using Ceva's theorem. Forum Geometricorum, vol. 16 (2016), 185-186.\\

Dr. Martin Celli.\\
Depto. de Matem\'aticas, Universidad Aut\'onoma
Metropolitana-Iztapalapa. M\'exico, D.F.\\
E-mail: cell@xanum.uam.mx

\end{document}